\documentclass[11pt]{amsart}
\usepackage{mathrsfs}
\usepackage{threeparttable}
\usepackage{amsfonts}
\usepackage{amsmath}
\usepackage{amssymb}
\usepackage{enumerate}
\usepackage{hyperref}
\usepackage{latexsym}
\usepackage{color}
\usepackage{enumitem}

\usepackage{enumerate}

\def\qed{\hfill $\Box$} \def\demo{\noindent{\it Proof.}~}

\newcommand\ZZ{\mathrm C} \newcommand\D{\mathrm{D}}

\newcommand\Cay{\mathrm{Cay}}

\newcommand\N{\mathbf{N}} \newcommand\Z{\mathbf{Z}}

     \newcommand\PSL{\mathrm{PSL}}

\newcommand\M{\mathrm{M}} \newcommand\Q{\mathrm{Q}}

\newtheorem{theorem}{Theorem}[section]

\newtheorem{lemma}[theorem]{Lemma}

\newtheorem{proposition}[theorem]{Proposition}
\theoremstyle{definition}

\def \Inv{{\rm{Inv}}}
\def\o{{\rm{o}}} \def\Syl{\hbox{\rm Syl}}
\def\lg{\langle} \def\rg{\rangle} 
\begin{document}
\title[Subgroup perfect codes]{Subgroup Perfect Codes of $\mathcal{A}_t$-Groups and Their Applications}
\thanks{Corresponding author: Bin Bin Li}
\thanks{This work was supported by  NNSFC (12301446),  NSF of Guangxi (2025GXNSFAA069013) and the Special Foundation for Guangxi Ba Gui Scholars}
\author[H.Y. Chen, B.B. Li, J.J. Li \and H. Yu]{%
Hu Ye Chen, Bin Bin Li, Jing Jian Li and Hao Yu}
\address{}
\address{School of Mathematics and Information Science \& Guangxi Base \\
   Tianyuan Mathematical Center in Southwest China  \\
\& Guangxi Center for Mathematical Research   \\
\& Center for Applied Mathematics of Guangxi(Guangxi University) \\
Guangxi University \\
Nanning, Guangxi 530004, P. R. China.}
\email{chenhy280@gxu.edu.cn (H.Y. Chen); \allowbreak libb@st.gxu.edu.cn (B.B. Li); \allowbreak lijjhx@gxu.e\allowbreak{du.cn} (J.J. Li); \allowbreak  haoyu@gxu.edu.cn (H. Yu).}

\begin{abstract}
A subset $C$ of the vertex set of a graph $\Gamma$ is called a perfect code in $\Gamma$ if every vertex of $\Gamma$ is at distance no more than 1 to exactly one vertex of $C$. A subgroup $H$ of a group $G$ is called a subgroup perfect code of $G$ if $H$ is a perfect code in some Cayley graph of $G$. Recently, Zhang reveals that the study of subgroup perfect codes of finite groups naturally reduces to the case of $p$-groups, especially $2$-groups.
Based on the combined works of Berkovich, Janko and Zhang, every $p$-group is an $\mathcal{A}_t$-group.
In this work, we establish a complete classification of subgroup perfect codes of  $\mathcal{A}_t$-groups for $t \in\{0, 1\}$.
Moreover, subgroup perfect codes of finite groups with abelian Sylow $2$-subgroups are also characterized.

\vskip 5pt

\end{abstract}
\maketitle

\parskip 5pt

\section{Introduction}
In this work, all groups considered are finite, and all graphs considered are finite, simple and undirected.
Given a group $G$ and an inverse-closed subset $S$ of $G \setminus \{1\}$, the \textit{Cayley graph} $\Cay(G, S)$ is the graph with vertex set $G$, where two distinct vertices $x$ and $y$ are adjacent if and only if $yx^{-1} \in S$.
A subset $C$ of the vertex set of graph $\Gamma$ is called a \textit{perfect code} \cite{j} in $\Gamma$ if every vertex of $\Gamma$ is at distance no more than 1 to exactly one vertex of $C$ (in particular, $C$ is an independent set of $\Gamma$).
Perfect codes are equivalently known as efficient dominating sets \cite{dej} or independent perfect dominating sets \cite{l}.
The study of perfect codes, particularly their realization within Cayley graphs, is an active area of research. Background and foundational results can be found in \cite[Section 1]{h}; for more recent developments, see \cite{den, f, ku}.
When a perfect code $C$ in a Cayley graph $\Cay(G, S)$ additionally forms a subgroup of $G$, it possesses both combinatorial properties from the graph and algebraic structure from the group.
This interaction leads to the problem of characterizing which subgroups of $G$ can be realized as perfect codes in some Cayley graph.
Huang, Xia, and Zhou \cite{h} initiated the systematic study of this problem by introducing the notion of subgroup perfect codes.

A subgroup $H$ of $G$ is called a \textit{subgroup perfect code} if there exists some inverse-closed subset $S$ of $G\setminus\{1\}$ such that $H$ is a perfect code in Cayley graph $\Cay(G, S)$. Clearly, the trivial subgroups $1$ and $G$ are perfect codes in the Cayley graphs $\Cay(G, G \setminus \{1\})$ and $\Cay(G, \emptyset)$, respectively.
Given a normal subgroup $H$ of $G$, Huang et al.\cite{h} provided a necessary and sufficient condition for $H$ to be a subgroup perfect code of $G$.
This work was extended by Chen, Wang and Xia \cite{ch}, who established a series of equivalent conditions for a subgroup $H$ of $G$ to be a subgroup perfect code of $G$ (see Proposition \ref{chen}).
Their work also revealed that every finite group $G$, except for cyclic $2$-groups and generalized quaternion $2$-groups, admits a nontrivial subgroup as its subgroup perfect code.
A group $G$ is said to be \textit{code-perfect} if every subgroup of $G$ is a subgroup perfect code \cite{m}. Ma, Walls, Wang and Zhou \cite{m} proved that a group is code-perfect if and only if it contains no elements of order 4.
Consequently, every finite odd-order group is automatically code-perfect.
This fact therefore redirects research attention to the core problem: characterizing subgroup perfect codes of groups with even-order. Recently, Zhang \cite{z} showed that this problem can be reduced to the study of subgroup perfect codes of $2$-groups(see Proposition \ref{zhang}). Thus, studying subgroup perfect codes of $p$-groups is particularly important.

However, classifying subgroup perfect codes of $p$-groups remains extremely difficult, primarily since the count of non-isomorphic $p$-groups of order $p^n$ grows very rapidly; see the asymptotic formula below, which was established by Higman and Sims \cite{Higman1,Higman2,Sims}:
\[
f(n,p) = p^{n^3(2/27 + O(n^{-1/3}))} \quad \text{(as $n \rightarrow \infty$)}.
\]
It is well-known that every $p$-group is an $\mathcal{A}_t$-group for $t\ge 0$, where the theory of $\mathcal{A}_t$-groups was introduced by Berkovich and Janko \cite{Berkovich-Janko-book}, and further developed by Zhang \cite{ZQH}.
For a positive integer $t$, a $p$-group $G$ is called an \textit{$\mathcal{A}_t$-group} if it contains a nonabelian subgroup of index $p^{t-1}$, but all its subgroups of index $p^t$ are abelian; for $t=0$, an $\mathcal{A}_0$-group \cite{ZQH} is defined to be an abelian $p$-group.
Specifically, $\mathcal{A}_1$-groups are the well-known \textit{minimal nonabelian $p$-groups} (where all proper subgroups are abelian, but the group itself is non-abelian). This suggests that studying subgroup perfect codes of $p$-groups is equivalent to studying subgroup perfect codes of $\mathcal{A}_t$-groups.

In this work, we characterize subgroup perfect codes of $\mathcal{A}_t$-group where $t \in \{0, 1\}$. Our main theorems are as follows.

\begin{theorem}\label{main}
Let $G$ be a $p$-group and $H$ a subgroup of $G$, where $p$ is a prime.
If $G$ is an $\mathcal{A}_t$-group where $t \in \{0, 1\}$,
then $H$ is a subgroup perfect code of $G$ if and only if either $H\in\{1, G\}$ or one of the following holds:
\begin{enumerate}[font=\normalfont]
\item If $p$ is an odd prime, then $G$ is code-perfect;

\item If $(t,p)=(0,2)$, then $H \cap \Phi(G) \leqslant \Phi(H)$;

\item If $(t,p)=(1,2)$, then either
       \begin{enumerate}[label=\upshape(\alph*), align=left,  widest=(b)]
      \item $H=\lg x \rg$, where $x$ is a nonsquare element of $G$ with $G\not\cong \Q_8$; or
      \item $H$ is noncyclic, and $(H,G)\in \{(H_0,\D_8), (H_1,G_1), (H_2,G_2)\}$ where $H_0\cong \ZZ_2\times \ZZ_2$,
            $G_1=\lg a, b \mid a^2 = b^{2^m} = c^2 = 1,\ [a, b] = c,\ [a,c]=[b,c]=1 \rg$,
            $G_2=\lg a, b \mid a^{2^n} = b^{2^m} = c^2 = 1,\ [a, b] = c,\ [a,c]=[b,c]=1 \rg$,
            $H_1\in \{ \lg a c^s, b^2\rg,\ \lg a b^{2j} c^s, b^{2^k r} c \rg, \allowbreak \lg a b^t, c \rg,\ \lg a^t b^d, c \rg\}$,
            $H_2 \in  \{ \langle a^d c^s, b^2 \rangle, \langle a^d b^{2j} c^s, b^{2^k r} c \rangle, \allowbreak \langle a^t b^d c^s, a^2\rangle, \langle a^t b^d c^s, a^{2^l r} c \rangle, \lg a^d b^t, c \rg$,\  $\lg a^t b^d, c \rg\}$, with  $2\le n \le m$, $t\ge 0$ and $j\ge 0$, $d$ and $r$ are odd integers, $1\le k < m$, $2^k\mid 2^{n}j$, $1\le l < n$ and $s\in \{0,1\}$.
      \end{enumerate}
\end{enumerate}
\end{theorem}

\begin{theorem}\label{abelian}
Let $G$ be a finite group with a nontrivial abelian Sylow $2$-subgroup and $H\leqslant G$.
Let $Q\in\Syl_2(H)$ and $P\in\Syl_2(G)$ such that $Q\leq P$.
Then $H$ is a subgroup perfect code of $G$ if and only if $Q \cap \Phi(P) \leqslant \Phi(Q)$.
Moreover,
\begin{itemize}[font=\normalfont]
\item [(1)] If  $G$ is simple, then $G$ is code-perfect; and

\item [(2)] If $G$ is a minimal nonabelian group and $P \neq G$, then $G$ is code-perfect if and only if either $P\lhd G$ or $P\cong \ZZ_2$; furthermore, if $P\ntriangleleft G$, then $H$ is a subgroup perfect code of $G$ if and only if $Q\in\{1,P\}$.
\end{itemize}
\end{theorem}

After this introductory section, some notations, basic definitions and useful facts will
be given in Section \ref{sec2}, and Theorem \ref{main} and \ref{abelian} will be proved in Section \ref{sec3}.

\section{Preliminaries}\label{sec2}

Notations and terminologies used in the paper are standard and can be found in \cite{hu}.
For example, we use $\Phi(G)$ to denote the \emph{Frattini subgroup} of group $G$, $\Syl_p(G)$ to denote the set of all Sylow $p$-subgroup of group $G$ with a prime $p$, and $d(G)$ to denote the minimal number of generators for group $G$.
Moreover, for an element $x \in G$, the order of $x$ is written as $\o(x)$, and $x$ is called a \textit{square} if there exists an other element $y \in G$ such that $x = y^2$.
For a 2-group $G$, set $\Inv(G) = \{ x \in G \mid x^2 = 1 \}$.
Then $\Omega_1(G)=\langle\Inv(G)\rangle$.

About the Frattini subgroup of $G$, the following result is useful.
Remind that $G^n=\langle x^n \mid x \in G \rangle$ for some integer $n\ge0$.

\begin{proposition}\cite[III, Satz 3.14]{hu}\label{G^2}
Let $G$ be a $p$-group with a prime $p$.
Then $\Phi(G)=G' G^p$ and $G/\Phi(G)$ is an elementary abelian $p$-group.
In particular, if $p=2$, then $\Phi(G)=G^2$.
\end{proposition}

\begin{proposition}\cite[Theorem 1.1]{m}\label{odd}
A group is code-perfect if and only if it has no elements of order $4$.
\end{proposition}

\begin{proposition}\cite[Theorem 1.2]{z} \label{zhang}
Let $G$ be a finite group and $H\leq G$.
Set $Q\in\Syl_2(H)$ and $P\in\Syl_2(\N_G(Q))$.
Then $H$ is a subgroup perfect code of $G$ if and only if $Q$ is a subgroup perfect code of $P$.
\end{proposition}

\begin{proposition}\cite[Theorem 2.11]{h} \label{dihedral}
Let $G=\langle a, b \mid a^n=b^2=1, (ab)^2=1 \rangle\cong \D_{2n}$. Then
the subgroup $H$ of $G$ is a subgroup perfect code of $G$ if and only if either $H\leq \langle a \rangle$ and $|H|$ or $n/|H|$ is odd; or $H\nleq \langle a \rangle$.
\end{proposition}

\begin{proposition}\cite[Theorem 1.2]{ch} \label{chen}
Let $G$ be a group and $H\leqslant G$. Then the following statements are equivalent:
\begin{itemize}[font=\normalfont]
\item [(1)]  $H$ is a subgroup perfect code of $G$;

\item [(2)]  there exists an inverse-closed right transversal of $H$ in $G$;

\item [(3)]  for each $x \in G$ such that $x^2\in  H$ and $|H|/|H \cap H^x|$ is odd, there exists $y \in Hx$ such that $y^2 = 1$;

\item [(4)]  for each $x\in G$ such that $HxH = Hx^{-1}H$ and $|H|/|H \cap H^x|$ is odd, there exists $y \in Hx$ such that $y^2 = 1$.
\end{itemize}
\end{proposition}

\vskip 3mm
Based on Proposition \ref{chen}, the following lemma is useful.

\begin{lemma} \label{cyclic}\label{omega}
Let $G$ be a group and $H$ a nontrivial subgroup of $G$.
Then we have:
\begin{itemize}[font=\normalfont]
\item [(1)]  if $H=\langle g^2\rangle$ is a $2$-group where $1\neq g\in G$, then $H$ is not a subgroup perfect code of $G$;

\item [(2)]  if $G$ is a $2$-group and $\Inv(G)\subseteq H$, then $H$ is not a subgroup perfect code of $G$.
\end{itemize}
\end{lemma}

\proof (1)
Suppose that $H=\langle g^2 \rangle$, where $\o(g)=2^k$ for $k\geqslant 2$.
Then $\langle g \rangle= H \cup Hg$. Since there exists only one involution in $\lg g \rg$ and also in $H$, we get that $Hg$ contains no involution.
 By Proposition \ref{chen}.(3), $H$ is not a subgroup perfect code of $G$.

(2) Suppose that $G$ is a 2-group and $\Inv(G)\subseteq H$. Then $H< \N_G(H)$ as $G$ is nilpotent.
There exists an element $x\in \N_G(H)\setminus H$ such that $Hx^2=H$, i.e., $x^2\in H$.
Since $H^x\cap H=H$ and the coset $Hx$ contains no involution, $H$ is not a subgroup perfect code of $G$, by Proposition \ref{chen}.(3).
\qed

\vskip 3mm
Miller and Moreno\cite{Miller} revealed the structure of minimal non-abelian groups.

\begin{proposition}\cite{Miller}\label{Miller}
Let $G$ be a minimal non-abelian group. Then precisely one of the following holds:
\begin{itemize}[font=\normalfont]
\item [(1)] $G$ is an $\mathcal{A}_1$-group (i.e., minimal nonabelian $p$-group with a prime $p$);
\item [(2)] $G = P {:} Q$ with $P$ an elementary abelian Sylow $p$-subgroup, $Q$ a cyclic Sylow $q$-subgroup, and $p \neq q$ are primes.
\end{itemize}
\end{proposition}

\begin{proposition}\cite[Lemma 2.3]{A2}\label{pro1}
Let $G$ be a $2$-group. Then $G$ is minimal nonabelian if and only if $d(G)=2$ and $|G'|=2$; if and only if $d(G)=2$ and $\Z(G) =\Phi(G)$.
\end{proposition}

The systematic study of $\mathcal{A}_t$-groups was initiated by R\'edei \cite{Redei}, who established the foundational classification for $t=1$.

\begin{proposition}\cite{Redei} \label{A1}
Let $G$ be an $\mathcal{A}_1$-group and a $2$-group.
Then $G$ is one of the following:
\begin{itemize}[font=\normalfont]
\item [(1)] the quaternion group $\Q_8$;

\item [(2)] the metacyclic group $\M_2(n_1,m_1)=\langle a,b \mid a^{2^{n_1}}=b^{2^{m_1}}=1, b^{-1}ab=a^{1+2^{n_1-1}}  \rangle=\langle a \rangle{:}\langle b \rangle$ with $n_1\geq 2$; or

\item [(3)] the non-metacyclic group $\M_2(n_2,m_2,1)=\langle a,b \mid a^{2^{n_2}}=b^{2^{m_2}}=1, [a,b]=c, [a,c]=[b,c]=1   \rangle$ with order $2^{n_2+m_2+1}$, where $n_2+ m_2\ge3$.
\end{itemize}
\end{proposition}

\vskip 3mm

Using Proposition \ref{A1}, we get the following lemma.

\begin{lemma}\label{pro2}
With the notation of  Proposition \ref{A1}, if $|G|>8$, then we have either
$G\cong \M(n_1,m_1)$, where $n_1\geq 2$ and $m_1+n_1\ge 4$; or $G\cong \M(n_2,m_2,1)$, where $m_2+n_2\ge3$.
Moreover,
\begin{itemize}[font=\normalfont]
\item [(1)] if $G\cong \M(n_1,m_1)$, then $\Inv(G)=\Omega_1(G)\cong \ZZ_2 \times \ZZ_2$;

\item [(2)] if $G\cong \M(n_2,m_2,1)$, then $\Inv(G)=\Omega_1(G)\cong \ZZ_2^3$.
\end{itemize}
\end{lemma}
\proof Under the hypothesis, Proposition \ref{A1} implies that either $G\cong\M(n_1,m_1)$ or $G\cong\M(n_2,m_2,1)$ where $n_1\ge 2$, $n_1+m_1\ge 4$ and $m_2+n_2\geq 3$. Hence, the proof is divided into the following two cases.
\vskip 3mm
{\it Case 1: $G\cong \M(n_1,m_1)$.}
\vskip 3mm

By Proposition \ref{A1}, $G=\lg a,b \mid a^{2^{n_1}}=b^{2^{m_1}}=1, b^{-1}ab=a^{1+2^{n_1-1}} \rg =\lg a \rg {:} \lg b \rg$, where $n_1\ge2$ and $n_1+m_1\geq 4$.
Set $g=a^i b^j$, an involution of $G$, with two integers $i,j$.
Set $g=a^i b^j$, an involution of $G$, with two integers $i,j$.
Then $$g^2=(a^i b^j)^2=a^i(a^i)^{b^{2^{m_1}-j}} b^{2j}=a^{i+i(1+2^{{n_1}-1})^{2^{m_1}-j}}b^{2j}=1,$$
which implies  $i+i(1+2^{n_1-1})^{2^{m_1}-j}\equiv0\pmod{2^{n_1}}$ and $2j\equiv0\pmod{2^{m_1}}$.
If $2^{m_1}-j$ is even, then $i+i(1+2^{n_1-1})^{2^{m_1}-j}\equiv 2i\equiv0\pmod{2^{n_1}}$ and so $g\in \langle a^{2^{n_1-1}}, b^{2^{m_1-1}}\rangle=\langle a^{2^{n_1-1}}\rangle \times \langle b^{2^{m_1-1}}\rangle$.
Now assume that $2^{m_1}-j$ is odd.
By the equation $2j\equiv0\pmod{2^{m_1}}$, we get $m_1=1$ and $j=1$, which implies $b=b^{2^{m_1-1}}$ and $n_1\ge3$.
Then $2i+i2^{n_1-1}\equiv0\pmod{2^{n_1}}$, which implies $i(1+2^{n_1-2})\equiv0\pmod{2^{n_1-1}}$.
Since $n_1\ge3$, we get that $2^{n_1-2}$ is even and so $1+2^{n_1-2}$ is odd, which implies $i\equiv0\pmod{2^{n_1-1}}$.
Then $g\in \langle a^{2^{n_1-1}}\rangle \times \langle b^{2^{m_1-1}}\rangle$.
We therefore establish
$\Omega_1(G)=\langle a^{2^{n_1-1}}\rangle \times \langle b^{2^{m_1-1}} \rangle \cong \ZZ_2 \times \ZZ_2$.

\vskip 3mm
{\it Case 2: $G\cong \M(n_2,m_2,1)$.}
\vskip 3mm

By Proposition \ref{A1} again, $G=\langle a,b \mid a^{2^{n_2}}=b^{2^{m_2}}=1, [a,b]=c, [a,c]=[b,c]=1  \rangle$, where $m_2+n_2\ge3$.
By Propositions \ref{G^2} and \ref{pro1}, we conclude that $G^2=\Phi(G)=\Z(G)$ and $c\in \Z(G)$ is an involution.
Since $[a,b]=c$, the element $a$ normalizes the subgroup $\langle b\rangle \times \langle c\rangle$.
Then $G=\langle a\rangle (\langle b\rangle \times \langle c\rangle)$.
Combining with the fact $|G|=2^{m_2+n_2+1}$, we obtain $\langle a\rangle \cap (\langle b\rangle \times \langle c \rangle)=1$.
Set $g=a^i b^j c^k$, an involution of $G$, with two integers $i,j$ and $k\in\mathbb{Z}_2$.
Then the relation $ab=bac$ implies
$$g^2 = (a^i b^j c^k)^2=(a^i b^j)^2 =a^{2i}b^{2j}c^{ij}= 1.$$
If both \(i\) and \(j\) are odd, then \(\langle a \rangle \cap \langle b,c \rangle = 1\) implies \(a^{-2i} = b^{2j}c = 1\), contradicting \(\langle b \rangle \cap \langle c \rangle = 1\). Thus \(i\) or \(j\) is even, so \(a^{2i}b^{2j} = 1\), giving \(2i \equiv 0 \pmod{2^{n_2}}\) and \(2j \equiv 0 \pmod{2^{m_2}}\).
Then $g\in \langle a^{2^{n_2-1}}\rangle \times \langle b^{2^{m_2-1}}\rangle \times \langle c\rangle$.
We therefore establish
$\Omega_1(G)=\langle a^{2^{n_2-1}}\rangle \times \langle b^{2^{m_2-1}} \rangle \times \langle c\rangle \cong \ZZ_2^3$. \qed
\vskip 3mm
The characterization of finite groups with abelian Sylow $2$-subgroups was established by Walter \cite{Walter}.
Combing with \cite{hu,JK,Ree}, we provide a complete characterization of finite nonablian simple groups with abelian Sylow $2$-subgroups.
\begin{proposition}\cite{Walter,hu,JK,Ree} \label{simple}
Let $G$ be a finite nonablian simple group with an abelian Sylow $2$-subgroup $P$.
Then one of the following holds true:
\begin{itemize}[font=\normalfont]
\item [(1)] $G \cong \PSL(2,2^n)$ and $P \cong \ZZ_2^n$ {\rm(}see \cite[II Satz 8.10 \text{(a)}]{hu}{\rm )}, where $n > 1$;
\item [(2)] $G \cong \PSL(2,q)$, where $q \equiv \pm3\pmod{8}$, and $P \cong \ZZ_2^2$ {\rm(}see \cite[II Satz 8.10 \text{(b)}]{hu}{\rm )};
\item [(3)] $G$ is the Janko group $\mathrm{J}_1$,  and $P \cong \ZZ_2^3$ {\rm(}see  \cite{JK}{\rm)};
\item [(4)] $G$ is the Ree group  $\mathrm{^2G}_2(q)$, where $q = 3^{2n+1}$ and $n \geqslant 1$, and $P$ is elementary abelian {\rm(}see \cite[Theorem 8.3]{Ree}{\rm)}.
\end{itemize}
\end{proposition}
\section{Subgroup perfect codes of $\mathcal{A}_t$-groups where $t=0$ or $1$.}\label{sec3}
In this section, we shall establish a complete characterization of subgroup perfect codes of $\mathcal{A}_t$-groups for $t\in \{0,1\}$. Let $G$ be an $\mathcal{A}_t$-group which is a $p$-group, and let $1 < H < G$, where $p$ is prime.
\begin{lemma}\label{t=0}
Suppose that $t=0$.
Then if $p$ is odd, then $G$ is code-perfect; if $p=2$, then $H$ is a subgroup perfect code of $G$ if and only if $H\cap \Phi(G)\leq \Phi(H)$.
\end{lemma}
\demo
Suppose that $p$ is an odd prime.
Then there exists no element of order 4 in $G$.
By Proposition \ref{odd}, $G$ is code-perfect, as desired.
So in what follows, we assume that $p=2$.
Now, $G$ is an abelian $2$-group, which implies that $G^2=\{g^2 \mid g\in G\}$.
By Proposition \ref{G^2}, we get that $\Phi(G)=\{g^2 \mid g\in G\}$.

Suppose that $H \cap \Phi(G)\nleq \Phi(H)$. Then there exists an element $g^2\in (H\cap \Phi(G))\setminus \Phi(H)$. Since $g^2 \notin \Phi(H)=\{h^2\mid h\in H\}$, we get that $g\notin H$ and $g^2\in H$.
Since $H \unlhd G$, we get $H=H^g$ and so $|H|=|H\cap H^g|$. We claim there exists no involution $x\in Hg$. Indeed, otherwise, $x=h_1 g$, where $h_1\in H$ and so $x^2 = (h_1 g)^2=h_1^2 g^2 =1$. Thus, $g^2=(h_1^2)^{-1} \in \Phi(H)$, a contradiction. By Proposition \ref{chen}.(3), $H$ is not a subgroup perfect code of $G$.

Suppose that $H\cap \Phi(G)\leq \Phi(H)$. Since $G$ is abelian, $H\unlhd G$. Let $g\in G$ such that $g^2\in H$. Then $g^2\in H\cap \Phi(G)\leq \Phi(H)$. Consequently, there exists $h\in H$ such that $g^2= h^2$ as $\Phi(H)=\{ h^2 \mid h\in H \}$. Since $(h^{-1}g)^2=h^{-2}g^2=1$, the element $h^{-1}g \in Hg$ satisfies $(h^{-1}g)^2=1$. By Proposition \ref{chen}.(3), $H$ is a subgroup perfect code of $G$. \qed

\begin{lemma}\label{t=1}
Suppose $t=1$.
Then if $p$ is odd, then $G$ is code-perfect;  if $p=2$, then $H$ is a subgroup perfect code of $G$ if and only if one of the following statements holds:
\begin{enumerate}[font=\normalfont]

\item $H=\lg x \rg$, where $x$ is a nonsquare element of $G$, except for $G \cong \mathrm{Q}_8$;

\item $H$ is noncyclic, and either $G\cong \D_8$;  or\\
$G=\lg a, b \mid a^{2^{n_2}} = b^{2^{m_2}} = c^2 = 1,\ [a, b] = c,\ [a,c]=[b,c]=1 \rg$ with $1\le n_2 \le m_2$, $n_2+ m_2\ge 3$, and
$H \in \{ \lg a c^s, b^2\rg, \lg a b^{2j} c^s, b^{2^k r} c \rg, \lg a b^t, c \rg, \lg a^t b^d, c \rg\}$ if  $n_2 = 1$, or
$H \in \{ \lg a^d c^s, b^2 \rg, \lg a^d b^{2j} c^s, b^{2^k r} c \rg, \lg a^t b^d c^s, a^2\rg, \lg a^t b^d c^s, a^{2^l r} c \rg, \lg a^d b^t, c \rg, \lg a^t b^d, c \rg  \}$ if $n_2\ge 2$,
 where $t$ and $j$ are nonnegative integers, $d$ and $r$ are odd integers, $1\le k < m_2$, $2^k\mid 2^{n_2}j$, $1\le l <n_2$ and $s\in \{0,1\}$.
\end{enumerate}
\end{lemma}
\demo \
Suppose that $p$ is an odd prime.
Then there exists no element of order 4 in $G$.
By Proposition \ref{odd}, $G$ is code-perfect, as desired.
So in what follows, we assume that $p=2$.
By Proposition \ref{A1}, we get $G\in \{\Q_8, \M_2(n_1,m_1), \M_2(n_2,m_2,1)\}$, where $n_1\ge 2$, $1\le n_2 \le m_2$ and $n_2 + m_2 \ge 3$. So in what follows, we divide it into three cases.
\vskip 3mm
{\it Case 1: $G\cong \Q_8$.}
\vskip 3mm
Suppose $G\cong \Q_8$. Since $G$ contains a unique involution, every nontrivial subgroup of $G$ contains $\Omega_1(G)$. Lemma \ref{omega}.(2) implies that $G$ admits only the trivial subgroup as a subgroup perfect code.

\vskip 3mm
{\it Case 2: $G \cong \M_2(n_1, m_1)$.}
\vskip 3mm
Now, $G=\langle a, b \mid a^{2^{n_1}}=b^{2^{m_1}}=1, b^{-1}ab=a^{1+2^{n_1-1}} \rangle$, where $n_1\geq 2$.
If $n_1=2$ and $m_1=1$, then $G\cong \D_8$, and the result follows from Proposition \ref{dihedral}.
So in what follows, we assume $n_1+m_1\geq 4$.
Then $|G|>8$. By Lemma \ref{pro2}.(1), we get $\Inv(G)=\Omega_1(G)\cong \ZZ_2 \times \ZZ_2$.

Suppose $H$ is noncyclic. Since $H$ is abelian, and $\Omega_1(H)\leq \Omega_1(G)\cong \ZZ_2\times \ZZ_2$, we get $4\leq |\Omega_1(H)|\leq |\Omega_1(G)|\leq 4$, which implies $\Omega_1(H)=\Omega_1(G)$. By Lemma \ref{omega}.(2), $H$ is not a subgroup perfect code of $G$.

Suppose that $H$ is cyclic. Set $H=\langle h \rangle$, where $h\in G$.
By Lemma \ref{omega}.(1), $H$ is a subgroup perfect code of $G$ only if $h$ is nonsquare element. So in what follows, we assume that $h$ is a nonsquare element of $G$.

Let $g\in \N_G(H)\setminus H$ such that $g^2\in H$. Then $K:=\langle H, g\rangle= H\langle g\rangle=H\cup Hg$ is noncyclic. If $K=G$, then $\Inv(K) \cap Hg \neq \emptyset$, which implies that there exists an involution in $Hg$. By Proposition \ref{chen}.(3), $H$ is a subgroup perfect code of $G$. Thus, suppose that $K<G$. Then $K$ is abelian as $G$ is an $\mathcal{A}_1$-group, which implies $\Inv(H)\subsetneq \Inv(K)$. Thus, $\Inv(K) \cap Hg \neq \emptyset$, which implies that there exists an involution in $Hg$. By Proposition \ref{chen}.(3), $H$ is a subgroup perfect code of $G$.

\vskip 3mm
{\it Case 3: $G \cong \M_2(n_2, m_2, 1)$.}
\vskip 3mm
Now, $G = \lg a, b \mid a^{2^{n_2}} = b^{2^{m_2}} = c^2 = 1,\, [a, b] = c,\,[a,c]=[b,c]=1 \rg$, where $1\le n_2\le m_2$ and $n_2 + m_2 \ge 3$. Then $|G|>8$. By Lemma \ref{pro2}.(2), we have $\Inv(G)=\Omega_1(G)\cong \ZZ_2^3$. Since $G$ is an $\mathcal{A}_1$-group, every proper subgroup of $G$ is abelian.

Suppose $H=\langle h \rangle$, where $h\in G$. Then by Lemma \ref{omega}.(1), we get that $H$ is a subgroup perfect code of $G$ only if $h$ is a nonsquare element. Further, with the same argument as the Case 2, we get that $H=\langle h \rangle$ is a subgroup perfect code of $G$ if and only if $h$ is a nonsquare element. So in what follows, we assume that $H$ is noncyclic and shall proof Lemma \ref{t=1}.(2).

We prove that if $H$ is a subgroup perfect code of $G$, then $\N_G(H)/H$ is cyclic and $d(H)=2$.  Assume that $H$ is a subgroup perfect code of $G$. Suppose $\N_G(H)/H$ is not cyclic. Then there exist distinct cosets $Hg_1, Hg_2 \in \Omega_1(\N_G(H)/H)$ with $g_1^2, g_2^2 \in H$. Since $H$ is a subgroup perfect code, Proposition \ref{chen}.(3) yields involutions $y_1 \in Hg_1$ and $y_2 \in Hg_2$. By Lemma \ref{pro2}, $\Omega_1(G) = \langle a^{2^{n-1}}\rangle \times \langle b^{2^{m-1}}\rangle \times \langle c\rangle$, implying that any product of two involutions is an involution. Since $\Omega_1(H) \cong \ZZ_2 \times \ZZ_2$, the set $H \cup Hy_1 \cup Hy_2$ contains exactly 11 involutions. However, this contradicts the fact that $|\Omega_1(G)| = 8$. Now suppose $d(H) = 3$. Since $H$ is abelian, $\Omega_1(H)\cong \ZZ_2^3$, forcing $\Omega_1(H)=\Omega_1(G)$. By Lemma \ref{omega}.(2), $H$ cannot be a subgroup perfect code of $G$, a contradiction. Thus both conditions hold.
\vskip 3mm
\textit{Subcase 3.1: $H\lhd G$.}
\vskip 3mm
Suppose that $H$ is a subgroup perfect code of $G$. Then $G=\N_G(H)$ and $G/H$ is cyclic, which implies $G'\leq H$.
Since $o(c)=2$ and $G/\lg c \rg$ is abelian, we get that $G'=\langle c \rangle$ and so $c\in H$.
Note that $ c $ is a nonsquare element of $G$. Then $ c \notin \Phi(H)=\{h^2\mid h\in H\} $. Consequently, $ c $ must serve as a generator of $ H $.
Given that $d(H)=2$, we may choose an element $x$ of $H$ such that $H=\langle x, c\rangle$. Since $o(c)=2$ and $H$ is abelian, we get $H=\langle x \rangle \times \langle c \rangle$. Next, we shall determine the form of the element $x$. Set $x=a^i b^j$ for integers $i,j$.
If $i,j$ both are even, then $H\leq \langle a^2, b^2, c \rangle =\Phi(G)$. Since $G/\Phi(G)$ is noncyclic, it follows that $G/H$ is noncyclic, a contradiction. Thus, at least one of $i$ or $j$ must be odd.

Next, we shall show that $H=\langle a^i b^j \rangle \times \langle c \rangle$ is a subgroup perfect code of $G$, where $i$ or $j$ is odd. Set an element $g=a^k b^l c^s\in G\setminus H$ such that $g^2\in H$, where $k,l,s$ are nonnegative integers.
Consider the group $K:=\langle H, g \rangle= H \langle g \rangle$. If $K=G$, then $G=H\cup Hg$, and $\Inv(K)=\Inv(G)\nsubseteq H$. Thus, $\Inv(K) \cap Hg \neq \emptyset$, which implies that there exists an involution in $Hg$. By Proposition \ref{chen}.(3), $H$ is a subgroup perfect code of $G$.
So in what follows, we assume that $K<G$. Then $K=\lg a^ib^j, c, g \rg$ is abelian and $\Phi(K)=\langle (a^i b^j)^2, g^2 \rangle \le H\cap \Phi(G)$.
Suppose that $d(K)=2$.
Then $K/\Phi(K)\cong \ZZ_2\times \ZZ_2$ and $H/\Phi(K)\cong \ZZ_2$ as $|K:H|=2$.
Note that $a^ib^j$ and $c$ are distinct nonsquare elements of $G$. Then $a^ib^j, c\notin \Phi(K)$.
However, since $a^ib^jc\notin \Phi(G)$, it follows that $a^ib^jc\notin \Phi(K)$, which implies that $H/\Phi(K)=\lg a^ib^j\Phi(K),\  c\Phi(K)\rg\cong \ZZ_2\times\ZZ_2$, a contradiction.
Thus, $d(K)=3$ and so $\Inv(K)=\Omega_1(K)\cong \ZZ_2^3$, which implies there exists an involution in $Hg$. By Proposition \ref{chen}.(3) again, $H$ is a subgroup perfect code of $G$.
\vskip 3mm
\textit{Subcase 3.2: $H\ntrianglelefteq G$.}
\vskip 3mm
Suppose that $H$ is a subgroup perfect code of $G$. Note that every maximal subgroup of $G$ is normal in $G$. Hence there exists a maximal subgroup $M$ of $G$ such that $H < M$ and $M=\N_G(H)$ as $M$ is abelian. Since $G / \langle c \rangle$ is abelian and $o(c) = 2$, it follows that $G' = \langle c \rangle$ and $c \in \Phi(G)$ by Proposition \ref{G^2}. Indeed, $c \notin H$. Otherwise, $c \in H$ would imply that $H / G'$ is normal in the abelian quotient $G / G'$, and hence $H \lhd G$, a contradiction. Given that $a^2, b^2 \in \Phi(G)$ and $G / \langle a^2, b^2, c \rangle \cong \ZZ_2^2$, we conclude that $\Phi(G) = \langle a^2, b^2, c\rangle$. As $G = \langle a, b\rangle$, this establishes that $M \in \{\langle \Phi(G), a \rangle, \langle \Phi(G), ab \rangle, \langle \Phi(G), b \rangle\}$. Clearly, $|M : \Phi(G)| = 2$.
By Proposition \ref{pro1},  $\Phi(G) = Z(G)$, which implies $H \not\leq \Phi(G)$ and $M = H \Phi(G)$.
Next, we consider two cases depending on the structure of $M$.

First, assume \(M = \langle \Phi(G), a \rangle=\lg a\rg\times\lg b^2\rg\times\lg c\rg\). Since \(M = \Phi(G) \cup \Phi(G)a\) and \(H \not\leq \Phi(G)\), there exists an element \(h \in H \cap(\Phi(G)a)\) of the form \(h = a^{i_0} b^{2j_0} c^{s_0}\), where \(i_0\) is odd, \(0 \leq j_0 < 2^{m_2-1}\), and \( s_0 \in\{0, 1\}\). Note that $H$ is abelian and $\Phi(H)\le \Phi(G)$.
Then $h\notin \Phi(H)$ and \(h\) serves as a generator of \(H\). Moreover, since \(H/H\cap\Phi(G) \cong M/\Phi(G)\cong \ZZ_2 \) and $h\notin H\cap \Phi(G)$, we have \(H = \langle H \cap \Phi(G), h \rangle\). Thus, we may choose an element \(x \in H \cap \Phi(G)\) of the form \(x = a^{2k_0} b^{2 l_0} c^{s_0'}\) such that \(H = \langle h, x \rangle\), where \(0 \leq k_0 < 2^{n_2-1}\), \(0 \leq l_0 < 2^{m_2-1}\) and \( s_0' \in\{0, 1\} \).
Then \(H = \langle h, h^{-2k_0'} x \rangle = \langle a^{i_0} b^{2j_0} c^{s_0}, b^{2l_0-4j_0 k_0'} c^{s_0'} \rangle\), where \(i_0 k_0' \equiv k_0 \pmod{2^{n_2-1}}\).
If \(l_0\) is even and $s_0'=0$, then $M/H=\langle b^2 H, c H \rangle \cong \ZZ_2\times \ZZ_2$, which contradicts $M/H$ is cyclic.
Thus, \(l_0\) is odd or \(s_0'=1\).
Since $a^{i_0} b^{2j_0} c^{s_0}\in H$ and
$(a^{i_0} b^{2j_0} c^{s_0})^{2^{n_2}} = a^{2^{n_2} i_0} b^{2^{n_2+1}j_0} c^{2^{n_2} s_0} = b^{2^{n_2+1} j_0},$
it follows that \(b^{2^{n_2+1} j_0} \in H\).
If \(s_0'=1\), then \(b^{2l_0-4j_0 k_0'} \notin H\) as \(c \notin H\), which implies $b^{2l_0-4j_0 k_0'}\notin \lg b^{2^{n_2+1} j_0}\rg$.
Thus, $\lg b^{2^{n_2+1} j_0}\rg< \lg b^{2l_0-4j_0 k_0'}\rg$, which implies \(2^k\mid 2^{n_2} j_0\) where $k$ is an integer such that $2^k\mid\mid (2l_0-4j_0 k_0')$.
To sum up, we get that: if \(M = \lg \Phi(G), a \rg=\lg a\rg\times\lg b^2\rg\times\lg c\rg\), then $H\in  \left\{ \lg a^d b^{2t} c^s, b^2 \rg, \lg a^d b^{2j} c^s, b^{2^k r} c \rg \right\}$, where $t$ and $j$ are nonnegative integers, $d$ and $r$ are odd integers, $1\le k < m_2$, $2^k\mid 2^{n_2}j$,  and $s\in \{0,1\}$.

Suppose $M = \langle \Phi(G), a^{\sigma}b \rangle$, where $\sigma\in \{0,1\}$. If $o(a)=2$, then $M=\langle a^{\sigma}b  \rangle \times \langle c \rangle$ and in this case, any proper noncyclic subgroup of $M$ is normal in $G$, which implies that $H\not< M = \langle a^{\sigma}b \rangle \times \langle c \rangle$, a contradiction. Thus, $o(a)>2$.
Since \(M = \Phi(G) \cup \Phi(G)a^{\sigma}b\) and \(H \not\leq \Phi(G)\), there exists an element \(h_1 \in H \cap (\Phi(G)a^{\sigma}b)\) of the form \(h_1 = a^{2i_1+\sigma} b^{j_1} c^{s_1}\), where \(j_1\) is odd, \(0 \leq i_1 < 2^{n_2-1}\), and \( s_1 \in \{0, 1\}\).
Then \(h_1\notin \Phi(H)\) as $\Phi(H)\le \Phi(G)$, which implies that $h_1$ serves as a generator of \(H\). Moreover, since \(H/H\cap\Phi(G) \cong M/\Phi(G)\cong \ZZ_2\), we have \(H = \langle H \cap \Phi(G), h_1 \rangle\). Then we may choose an element \(x_1 \in H \cap \Phi(G)\) of the form \(x_1 = a^{2 k_1} b^{2l_1} c^{s_1'}\) such that \(H = \langle h_1, x_1 \rangle\), where \(0 \leq k_1 < 2^{n_2-1}\), \(0 \leq l_1 < 2^{m_2-1}\), and \( s_1' \in\{0,1\}\).
Then \(H = \lg h_1, h_1^{-2 l_1'} x_1 \rg = \lg a^{2i_1+\sigma} b^{j_1} c^{s_1}, a^{2k_1-2(2 i_1+ \sigma) l_1'} c^{s_1'+\sigma l_1'} \rg\), where \(j_1 l_1' \equiv l_1 \pmod{2^{m_2-1}}\).
Since \(M/H\) is cyclic, this implies that \(k_1-(2 i_1+ \sigma) l_1'\) is odd or \(s_1'+\sigma l_1'\equiv 1 \pmod2\).
Moreover, if $s_1'+\sigma l_1'\equiv 1 \pmod2$, then since $c\notin H$, we get $a^{2k_1-2(2 i_1+ \sigma) l_1'}\neq 1$.
To sum up, we get that: if \(M = \lg \Phi(G), a^{\sigma}b \rg=\lg a^2 \rg\times\lg a^{\sigma}b\rg\times\lg c\rg\), then $n_2\geq 2$ and $H \in \left\{ \lg a^t b^d c^s, a^2\rg, \lg a^t b^d c^s, a^{2^l r} c \rg  \right\},$
where $t\geq 0$, $d$ and $r$ are odd integers, $1\le l < n_2$ and $s\in\{0,1\}$.

In conclusion,
\[
H \in \begin{cases}
\{ \langle a c^s, b^2\rangle, \langle a b^{2j} c^s, b^{2^k r} c \rangle \}, & \text{if } n_2 = 1; \\
 \{ \langle a^d c^s, b^2 \rangle, \langle a^d b^{2j} c^s, b^{2^k r} c \rangle, \langle a^t b^d c^s, a^2\rangle, \langle a^t b^d c^s, a^{2^l r} c \rangle  \}, & \text{if } n_2\ge 2,
\end{cases}
\]
where $t$ and $j$ are nonnegative integers, $d$ and $r$ are odd integers, $1\le k < m_2$, $2^k\mid 2^{n_2}j$, $1\le l < n_2$ and $s\in \{0,1\}$.

Let $\mathcal{H}_1=\left\{ \langle a c^s, b^2\rangle, \langle a b^{2j} c^s, b^{2^k r} c \rangle \right\}$ for $n_2=1$ and $\mathcal{H}_2= \langle a^d c^s, b^2 \rangle, \langle a^d b^{2j} c^s, b^{2^k r}$ $ c \rangle, \langle a^t b^d c^s, a^2\rangle, \langle a^t b^d c^s, a^{2^l r} c \rangle \} \text{ for } n_2\ge 2,$  where $t$ and $j$ are nonnegative integers, $d$ and $r$ are odd integers, $1\le k < m_2$, $2^k\mid 2^{n_2}j$, $1\le l < n_2$ and $s\in \{0,1\}$.

Suppose that $H\ntrianglelefteq G$ and $H\in \mathcal{H}_i$, where $i\in \{0,1\}$.
We shall show that \(H\) is a subgroup perfect code of \(G\).
Based on the preceding argument, we have that \(\N_G(H) = H\Phi(G)\) is a maximal subgroup of \(G\) and \(\N_G(H) \in \{\langle a, b^2,c\rangle, \langle a^2, b,c\rangle, \langle ab, a^2, c\rangle\}\).
Note that \(c \notin H\) since \(H \ntrianglelefteq G\).

Firstly, assume $n_2=1$ and $H\in \mathcal{H}_1$. If $H=\langle a b^{2j} c^s, b^{2^k r} c \rangle$, then $H\langle b^2\rangle=\N_G(H)$, which implies $\N_G(H)/H$ is cyclic. Note that $|\N_G(H):H|\ge2$.
Let \(g \in \N_G(H) \setminus H\) such that \(g^2 \in H\). Since the involution \(c \in \N_G(H) \setminus H\) and \(\N_G(H)/H\) is cyclic, we have \(Hg = Hc\). Therefore, \(Hg\) contains the involution \(c\). By Proposition \ref{chen}.(3), \(H\) is a subgroup perfect code of \(G\).
Observe that $\langle a c^s, b^2\rangle \langle c\rangle=\langle a, b^2, c\rangle$. By an same argument as above, we get that \(H\) is a subgroup perfect code of \(G\) for $H=\langle a c^s, b^2\rangle$.

Now suppose $n_2\ge 2$ and $H\in \mathcal{H}_2$. Assume \(H \in \{\langle a^d c^s, b^2 \rangle, \langle a^t b^d c^s, a^2 \rangle\}\). Then \(H\langle c\rangle = \N_G(H)\), and thus \(\N_G(H)/H\) is cyclic.
Note that $|\N_G(H):H|\ge2$.
Let \(g \in \N_G(H) \setminus H\) such that \(g^2 \in H\). Since the involution \(c \in \N_G(H) \setminus H\) and \(\N_G(H)/H\) is cyclic, we have \(Hg = Hc\). Therefore, \(Hg\) contains the involution \(c\). By Proposition \ref{chen}.(3), \(H\) is a subgroup perfect code of \(G\).

For \(H = \langle a^d b^{2j} c^s, b^{2^k r} c \rangle\)(resp. $H=\langle a^t b^d c^s, a^{2^l r} c \rangle$), we have \(H\langle b^2\rangle = \N_G(H)\)(resp. $H\langle a^2\rangle= \N_G(H)$), which implies \(\N_G(H)/H\) is cyclic. By an same argument as above, we get that \(H\) is a subgroup perfect code of \(G\). \qed

\vskip 3mm
\textbf{The proof of Theorem \ref{main}.} Let $G$ be an $\mathcal{A}_t$-group where $t\in\{0, 1\}$,  and let $H$ be a subgroup of $G$. If $H\in\{1,G\}$, then $H$ is a subgroup code of $G$; if $1<H<G$, then the result holds by Lemma \ref{t=0} for \(t = 0\) and Lemma \ref{t=1} for \(t = 1\).
\qed

\section{Applications}

\vskip 3mm
{\bf The proof of Theorem \ref{abelian}.}
Let $ G $ be a finite group with an nontrivial abelian Sylow $ 2 $-subgroup, and $ H \leq G $.
Fix a Sylow $ 2 $-subgroup $ Q $ of $ H $ and a Sylow $ 2 $-subgroup $ P $ of $ G $ such that $ Q \leq P $.
Observe that $ P \leq \N_G(Q) $. Since $ P $ is a Sylow $ 2 $-subgroup of $ G $, it follows that $ P $ is also a Sylow $ 2 $-subgroup of $ \N_G(Q) $. By Proposition~\ref{zhang}, $ H $ is a subgroup perfect code of $ G $ if and only if $ Q $ is a subgroup perfect code of $ P $. Furthermore, Theorem~\ref{main} implies that $ Q $ is a subgroup perfect code of $ P $ if and only if $Q \cap \Phi(P) \leq \Phi(Q)$. Thus, the first statement is valid.

Suppose that $G$ is simple. If $G$ is abelian, then $G\cong \ZZ_2$ and so $G$ is code-perfect. So in what follows, assume $G$ is nonabelian. By Proposition~\ref{simple}, every Sylow 2-subgroup of $G$ is elementary abelian.
Then $\Phi(P) = \Phi(Q) = 1$, which implies $Q \cap \Phi(P) = \Phi(Q)$.
By Theorem \ref{main}, we get that $Q$ is a subgroup perfect code of $P$. Moreover, applying Proposition \ref{zhang}, $H$ is a subgroup perfect code of $G$.
Then $G$ is code-perfect and Theorem \ref{abelian}.(1) holds.

\vskip 3mm

Suppose that $G$ is a minimal non-abelian group and $P\neq G$. Then $G$ is not a $p$-group.
By Proposition \ref{Miller}, $G$ is a semidirect product $G = Q_1{:}Q_2$, where $Q_1$ is an elementary abelian Sylow $q_1$-subgroup, $Q_2$ is a cyclic Sylow $q_2$-subgroup, $q_1\neq q_2$ are primes with $2\in\{q_1, q_2\}$.
If $q_1=2$, then $P=Q_1$ is elementary abelian, and by Proposition \ref{odd}, $G$ is code-perfect.
So in what follows, we assume $q_2=2$.  Then $P = Q_2$ is cyclic.
By Lemma \ref{omega}, a cyclic $2$-group $P$ admits only $1$ and $P$ as its subgroup perfect codes.
Moreover, applying Proposition \ref{zhang}, the subgroup $H$ of $G$ is a subgroup perfect code if and only if its Sylow $2$-subgroup $Q$ is either 1 or $P$. In particular, when \(P \cong \ZZ_2\), \(G\) has no elements of order \(4\). Then Proposition \ref{odd} implies that \(G\) is code-perfect, confirming Theorem \ref{abelian}.(2).
\qed

\noindent{\bf Declaration of competing interest}

The authors declare that they have no known competing financial interests or personal relationships that could have appeared to influence the work reported in this paper.

\end{document}